\documentclass[11pt]{amsart}
\usepackage{amsmath,amssymb}
\usepackage[curve,matrix,arrow]{xy}

\newtheorem{thm}{Theorem}[section]
\newtheorem{cor}[thm]{Corollary}
\newtheorem{lemma}[thm]{Lemma}
\newtheorem{prop}[thm]{Proposition}

\theoremstyle{definition}
\newtheorem{defn}[thm]{Definition}
\newtheorem{remark}[thm]{Remark}
\newtheorem{example}[thm]{Example}

\newcommand\wh[1]{\widehat{#1}}
\newcommand{\ovl}[1]{\overline{#1}}

\newcommand{\ls}[2]{{^{#1}\!{#2}}}

\newcommand{\qbox}[1]{\quad\hbox{#1}\quad}
\def\dst{\displaystyle}
\def\im{\operatorname{im}}

\def\id{\operatorname{Id}}

\def\End{\operatorname{End}\nolimits}
\def\Aut{\operatorname{Aut}\nolimits}
\def\Mat{\operatorname{Mat}\nolimits}

\def\GL{\operatorname{GL}\nolimits}

\def\spec{\operatorname{Spec}\nolimits}
\newcommand\prlim[1]{{\underset{{#1}}{\varprojlim}}}
\newcommand\inlim[1]{{\underset{{#1}}{\varinjlim}}}
\renewcommand\sup[1]{{\underset{{#1}}{\operatorname{sup}}}}
\renewcommand\liminf[1]{{\underset{{#1}}{\operatorname{lim inf}}}}

\def\Z{\mathbb Z}
\def\N{\mathbb N}
\def\Q{\mathbb Q}
\def\CN{\mathcal N}

\def\F{\mathbb F}
\def\R{\mathbb R}

\begin{document}

\title{The pro-$p$ group of upper unitriangular matrices}

\author{\sc Nadia Mazza}
\address
{Department of Mathematics and Statsitics\\ University of
  Lancaster\\ Lancaster \\LA1 4YF, UK} 
\email{n.mazza@lancaster.ac.uk}

\date{\today}


\begin{abstract}
We study the pro-$p$ group $G$ whose finite quotients
give the prototypical Sylow $p$-subgroup of the general linear groups
over a finite field of prime characteristic $p$. 
In this article, we extend the known results on the subgroup structure
of $G$. In particular, we give an explicit embedding of the Nottingham
group as a subgroup and show that it is selfnormalising. Holubowski
(\cite{holub1,holub0,holub2}) studies a free product $C_p*C_p$ as a
(discrete) subgroup of $G$ and we prove that its closure is
selfnormalising of infinite index in the subgroup of $2$-periodic
elements of $G$.
We also discuss change of rings: field
extensions and a variant for the $p$-adic integers, this latter linking
$G$ with some well known $p$-adic analytic groups.
Finally, we calculate the Hausdorff dimensions of some closed
subgroups of $G$ and show that the Hausdorff spectrum of $G$ is the
whole interval $[0,1]$ which is obtained by considering partition
subgroups only.
\end{abstract}

\maketitle

\smallskip


\smallskip


\section{Introduction}\label{sec:intro}

In this paper we investigate the pro-$p$ group whose finite quotients
give the prototypical Sylow $p$-subgroup of the general linear groups
over a finite field of prime characteristic $p$. For convenience, we
will consider an odd prime $p$ throughout the paper.

Sylow $p$-subgroups of finite general linear groups $\GL_n(q)$ for $q$ a
power of $p$ have been minutely analysed by Weir in the 50s
(\cite{weir}). His findings have
subsequently been exploited by many; in particular Bier
(\cite{bier1,bier2}) who extended some of Weir's results to the pro-$p$
group $G(q)$ of upper unitriangular matrices with coefficients in the
field $\F_q$ with $q$ elements. By {\em upper unitriangular}
matrix, we mean an upper triangular matrix with all diagonal
coefficients equal to $1$. In the first part of this article, we will
elaborate on Weir, Bier's and Holubowski's results
(\cite{holub1,holub0,holub2}), and we will focus on the subgroup
structure of $G(q)$, revisiting the notion of {\em partition subgroups}
considered by Weir. 
We will also discuss the embeddings
of a free product of the form $C_p*C_p$ as a discrete subgroup of
$G(p)$ and of the Nottingham group $\CN(q)$ (\cite{c1,c2,c3}) as a
closed pro-$p$ subgroup of $G(q)$. 
Then we will discuss how we can relate $G(q)$ and $G(p)$ for a field
extension $\F_q/\F_p$. In Section~\ref{sec:p-adic} we present a $p$-adic
version of the group $G(q)$ and briefly relate this group to some
well-known $p$-adic analytic groups (\cite{dix}).
In the last section of the paper we calculate the Hausdorff dimensions
(\cite{aber,bk,bs,falconer}) of the closed subgroups presented in the
preceding sections. 
A short appendix includes background about the automorphism groups of
the Sylow $p$-subgroups of the general linear groups $\GL_n(q)$ and
about the Hausdorff dimension for profinite groups.

\vspace{.6cm}
\begin{defn}\label{def1}
Let $q=p^f$ be a power of an odd prime number $p$, with $f\geq1$. 
For each $n\in\N$, let 
$G_n(q)$ be the Sylow $p$-subgroup of $\GL_n(q)$ formed by the upper
triangular matrices with diagonal coefficients all equal to $1$.
Let $V_n(q)=\F_q^n$ be the set of column vectors of size $n$ with
coefficients in $\F_q$.
\end{defn}

Note that for each $n>1$, we have
$$G_n(q)=\left(\begin{array}{c|c}
G_{n-1}(q)&V_{n-1}(q)\\
\hline
0_{1\times(n-1)}&1\end{array}\right)
\cong V_{n-1}(q)\rtimes G_{n-1}(q),$$
where $G_{n-1}(q)$ acts on $V_{n-1}(q)$ by left multiplication in the
obvious way.
Thus, the natural projections
\begin{equation}\label{eqn1}
\xymatrix{\dots\ar[r]&G_{n+1}(q)\ar[r]^{\pi_{n+1}}&
G_n(q)\ar[r]^{\pi_{n}}&\dots\ar[r]&G_2(q)\ar[r]^{\pi_{2}}&G_1(q)=1}
\end{equation}
form an inverse system. 
\begin{defn}\label{def:g}
Let $G(q)=\prlim nG_n(q)$ be the inverse limit of (\ref{eqn1}). The
group $G(q)$ is a pro-$p$ group, which we will call the {\em pro-$p$
  group of upper unitriangular matrices over $\F_q$}. 
If the prime power $q$ is clear from the context, we write simply $G$
and $G_n$ instead of $G(q)$ and $G_n(q)$ respectively.

For each $n\in\N$, let $\theta_n~:~G\to G_n$ be the universal map,
i.e. such that
$$\theta_m=\pi_{m+1}\cdots\pi_n\theta_n:G\to G_m
\qbox{for all $1\leq m<n$.}$$
Let also $N_n=\ker(\theta_n)$. Thus $N_n$ is the normal subgroup of $G$
formed by all the matrices whose upper left $n\times n$ diagonal block
is the identity matrix.
\end{defn}

From \cite[Section~1.2]{wilson}, a filter base for the topology on
$G$ is the set of open normal subgroups 
$$\mathcal B=\{\theta_n^{-1}(X)=XN_n~|~X\trianglelefteq G_n~,~n\in\N\}$$
while the set 
$\mathcal U=\{\theta_n^{-1}(X)=XN_n~|~X\subseteq G_n~,~n\in \N\}$ forms
a fundamental system of open neighbourhoods of the identity in $G$
(\cite[p.~26-27]{hig}). In particular, $G$ is countably based
(cf.~\cite[Proposition~4.1.3]{wilson}).  

For any $i,j\in\N$ let $e_{ij}$ denote the infinite elementary square
matrix whose unique nonzero coefficient is $(i,j)$ and is equal to $1$,
i.e. $(e_{ij})_{kl}=\delta_{ik}\delta_{jl}$ for all $i,j,k,l\in\N$. So,
if $a_1,\dots,a_f\in\F_q^\times$ generate $\F_q$ as $\F_p$-vector space, then 
$G=\langle1+a_ce_{ij}~|~1\leq i<j~,~1\leq c\leq f\rangle$, where we write
$1$ for the identity element of $G$ and $\F_q^\times$ for the
multiplicative group of nonzero elements of $\F_q$.
The set
$$\mathcal E=\{1+e_{i,i+1}~|~i\in\N\}$$
generates $G$ topologically and converges to $1$. Indeed, it generates a
dense subgroup of $G$, because 
$G_n=\langle 1+e_{i,i+1}~|~1\leq i<n\rangle$ for all $n\in\N$; moreover,
any open subgroup of $G$ contains all but a finite number of such
elements (\cite[Section 2.4]{RZ}). 

The metric on $G$ is defined as follows. 
Let $x,y\in G$ and fix a number $\epsilon\in(0,1)$, for instance 
$\epsilon=p^{-1}$. Then
\begin{equation}\label{eqn:dist}
d(x,y)=\epsilon^k\qbox{where}
k=\max\{n~|~y^{-1}x\in N_n\}.
\end{equation}
Then $d(~,~)$ is an ultrametric, i.e. subject to the axioms 
\begin{itemize}
\item $d(x,y)\geq0$ with equality if and only if $x=y$;
\item $d(x,y)=d(y,x)$;
\item $d(x,z)\leq\max\{d(x,y),d(y,z)\}$ for all $x,y,z\in G$\quad(the
  {\em ultrametric} axiom). 
\end{itemize}
For $x\in G$ and $k\geq0$, the open ball of centre $x$ and radius
$\epsilon^k$ is the set 
$B(x,\epsilon^k)=\{y\in G~|~d(x,y)<\epsilon^k\}$ of all the
elements $y$ of $G$ such that $y^{-1}x\in N_l$ for some integer
$l>k$. It differs from its closure only if $k$ is an integer, in
which case 
$\ovl{B(x,\epsilon^k)}=\{y\in G~|~y^{-1}x\in N_k\}=
B(x,\rho)=B(x,\epsilon^{k-1})$, 
for any $\epsilon^{k-1}\geq\rho>\epsilon^k$.
In particular,
$$N_n=\{y\in G~|~y\in N_n\}=\ovl{B(1,\epsilon^n)}=
B(1,\epsilon^{n-1})\qbox{for all $n\in \N$.}$$


\section{Partition subgroups of $G$}\label{sec:partition}\

Let $\gamma_n(G)$ denote the $n$-th term in the lower central
series of $G$, starting with $\gamma_1(G)=G$ and $\gamma_2(G)=[G,G]$,
the derived subgroup of $G$. The notation $[A,B]$ for subsets $A,B$ of a
group denotes the subgroup spanned by the commutators $[a,b]=a^{-1}a^b$
with $a\in A$ and $b\in B$, where we write
$$x^y=y^{-1}xy\qbox{and}\ls yx=yxy^{-1}.$$

At the basis of each computation, lays the ubiquitous commutator
relation
\begin{equation}\label{eqn2}
[1+ae_{ij},1+be_{kl}]=1+\delta_{jk}abe_{il}-\delta_{il}abe_{kj}
\qbox{for all $i,j,k,l\in\N$ and all $a,b\in\F_q$.}
\end{equation}

Extending work of Weir, Bier investigated a sub-family of the subgroups
that Weir called {\em partition subgroups} of $G$. The results she
proves are only for these partition subgroups, but it is easy to see
that they extend to all partition subgroups. For convenience, we
introduce the following definition. 

\begin{defn}
A {\em partition diagram} is a subset 
$$\mu=\{(r_i,c_i)\in\N^2~|~r_i<c_i\;,\;i\in\N\}$$
such that 
\begin{equation}\label{eqn:mu}
\qbox{whenever}(i,j)\in\mu\qbox{and}(j,k)\in\mu~,
\qbox{then}(i,k)\in\mu~.
\end{equation}
That is, a partition diagram is a collection of pairs of distinct positive
integers, which should be regarded as the coordinates of the nondiagonal
{\em squares} (or coefficients) in the matrices of $G$:  
$$(r,c)\in\mu\Longleftrightarrow1+e_{rc}\in G.$$
The corresponding {\em partition subgroup of $G$} is the subgroup
$$P_\mu=\{x\in G~|~x_{ij}=0\;,\;\forall\;(i,j)\notin\mu\}.$$
So, the constraint~(\ref{eqn:mu}) on the elements of $\mu$ reflects the
multiplication of the corresponding matrices (what Weir called
``completing the rectangle'') in $P_\mu$, namely 
$$(1+ae_{ij})(1+be_{jk})=1+ae_{ij}+be_{jk}+abe_{ik}.$$

If $\mu$ is such that for each $j\geq2$, if $(i,j)\in\mu$, then all
the pairs $(k,j)\in\mu$ for all $k\leq i$,
  then we call $\mu$ a {\em partition} and write it as
$$\mu=(\mu_2,\mu_3,\dots)\qbox{where}
\mu_j=\max\{l~|~(l,j)\in\mu\}\qbox{for all $j\geq2$.}$$
Then
$$P_\mu=\{x\in G~|~x_{ij}=0\;,\;\forall\;\mu_j<i<j\}$$
is formed by all the elements of $G$ whose $j$th column has
$(j-1-\mu_j)$ zeroes above the diagonal. 
A partition of the form 
$\mu=(0^{c-1},c,c+1,c+2,\dots)$ defines the partition subgroup $N_c$,
for any $c\geq1$ (and if $c=1$, then $N_1=G$). An exponent
``$\lambda^{s}$'' in a partition $\mu$ means $\lambda$ repeated $s$ times.

Given a partition diagram $\mu$, we denote $|\mu|$ 
its {\em shape}, i.e. the set of all squares on an infinite chessboard
$\N^2$ which consist of the possible nonzero squares $(i,j)$ in
$P_\mu$. 

A square $(i,j)$ {\em covers} $(k,l)$ if $(i,j)\neq(k,l)$ and if 
$k\leq i$ and $l\geq j$.
We say that $(i,j)$ {\em avoids} $|\mu|$ if $(i,j)$ covers a square
outside of $|\mu|$. 
\end{defn}

\begin{remark}
In \cite{bier1,bier2}, Bier only considers partitions. Moreover, she
takes the ``complementary'' definition of a partition than the one we
take here. That is, the parts in a partition denote the number of zeroes
above the diagonal. Instead, we have chosen to use the same convention
as Weir in \cite{weir}, in order to include the more general partition
subgroups defined by partition diagrams. 
\end{remark}

If $\mu\subseteq\N^2$ is a partition diagram, we call a {\em
  subpartition (diagram)} of $\mu$ a subset of $\mu$ which is a
partition (diagram) on its own. So a partition diagram is a lattice. That is
(\cite[Section~8.2]{hall2}), given any two subpartitions diagrams
$\mu_1$ and $\mu_2$ of $\mu$, their union and intersection are also
subpartition diagrams of $\mu$. The {\em union} of two partition
diagrams is the smallest partition diagram which contains them
(i.e. obtained by ``completing the rectangles'' in Weir's terminology),
whilst their {\em intersection} is their set intersection.
In particular, if
$\mu_1=(u_2,u_3,\dots)$ and $\mu_2=(v_2,v_3,\dots)$ are subpartitions of
$\mu$, then
\begin{align*}
\mu_1\cup\mu_2=(m_2,m_3,\dots)&\qbox{where}
m_j=\max\{u_j,v_j\}\qbox{and}\\
\mu_1\cap\mu_2=(n_2,n_3,\dots)&\qbox{where}
n_j=\min\{u_j,v_j\}.
\end{align*}
It follows that each partition diagram has a unique maximal subpartition
$$\mu_{\max}=\bigcup\mathfrak P_\mu
\qbox{where}\mathfrak P_\mu=\{\mu'\subseteq\mu~|~\mu'\;
\hbox{is a subpartition of $\mu$}~\}.$$
We say that a partition diagram $\mu$ {\em converges} to a partition if
there exists $n\geq2$ such that for any $(r,c)\in\mu$ with $c\geq n$,
then $(i,c)\in\mu$ for all $1\leq i\leq r$. That is, $\mu$ becomes a
partition for $n$ large enough. The {\em trivial partition} is the
partition $(0^{\aleph_0})$, where $\aleph_0$ is the cardinality of $\N$.

\cite[Theorem 2]{weir} describes the partition diagrams $\mu$ which
define normal subgroups $P_\mu$: namely $\mu$ is a partition and the
boundary of $|\mu|$ should move monotonically downward to the right.
The point is that if $P_\mu\trianglelefteq G$ and $(r,c)\in\mu$, then
conjugation by any $1+e_{ir}$ and $1+e_{cj}$ implies that $(i,c)$ and
$(r,j)$ must also be in $\mu$ for all $i\leq r$ and all $j\geq c$,
i.e. $\mu$ contains all the squares covered by $(r,c)$. So $\mu$ must be a
partition, and its ``boundary'', determined by all the squares $(i_c,c)$
with $i_c=\max\{r~|~(r,c)\in\mu\}$, must give an increasing
sequence $i_2\leq i_3\leq i_4\dots$.

A {\em rectangular} partition subgroup $P_\mu$ is a normal subgroup of
$G$ for $\mu$ of the form $\mu=(0^c,d^{\aleph_0})$,
where $0<d\leq c$ for some $c\in\N$ (we could extend to $d=0$ by
admitting the trivial subgroup of $G$ as a rectangular partition subgroup).
The shape of such $|\mu|$ explains the terminology. 
If $p>2$, then the maximal abelian (and characteristic) subgroups of $G$
have this form, with $d=c$ (\cite[Theorem 6]{weir}). That is, 
$$\mu=
(0^c,c^{\aleph_0})
\qbox{and}
P_\mu=\left(\begin{array}{r|r}I_{c+1}&\hbox{\Large*}\\
\hline0&I_\infty\end{array}\right)$$
where the coefficients in the block $\big(${\Large*}$\big)$ can take any value in
$\F_q$.

Extending Pavlov (\cite{pavlov}) and Weir's (\cite{weir})
results, Bier proves that the automorphism group of $G$ is generated by
three types of continuous automorphisms: inner, diagonal
(i.e. conjugation by an infinite diagonal matrix), and those induced by
field automorphisms. Furthermore, shifts are surjective group
homomorphism, where for $d\in\N$, the $d$th shift of $x\in G$ is the
matrix $x[d]$ obtained by deleting the first $d$ rows and columns of
$x$.

\begin{defn}\label{def:periodic}
We call a matrix $x\in G$ {\em periodic} (of period $d$) if there exists
$d\in\N$ such that $x=x[d]$. A subgroup $H\leq G$ is periodic (of period
$d$) if every element of $H$ is periodic (of period $d$).
\end{defn}

Here is a summary of Bier and Weir's results as they apply to $G=G(q)$. 

\begin{prop}\label{prop:bier-weir}\

\begin{enumerate}
\item Partition subgroups are closed.
\item A partition subgroup $P_\mu$ is open if and only if the
  partition diagram $\mu$ is such that there exists $N\in\N$ for which
  $(i,j)\in\mu$ for all $1\leq i<j$ and for all $j\geq N$.  
\item Let $H$ be a closed subgroup of $G$.
The following are equivalent.
\begin{enumerate}
\item $H$ is a normal subgroup of $G$.
\item $H$ is a normal partition subgroup of $G$.
\item $H$ is a characteristic subgroup of $G$.
\item $H$ is a partition subgroup defined by an increasing
  partition $\mu$, i.e. such that $\mu_j\leq\mu_{j+1}$. 
\end{enumerate}
If $H$ satisfies these equivalent conditions, we call $H$ a {\em normal
  partition} subgroup.
\item Given a partition diagram $\mu$, the normal core
$\cap_g\ls gP_\mu$ of $P_\mu$ is the partition subgroup 
$$\cap_g\ls gP_\mu=P_{\mu'}\qbox{where\;
$\mu'=(\mu'_2,\mu'_3,\dots)$\;with\;
$\mu'_j=\min\{i~|~(i,k)\in\mu~,~\forall~k\geq j\}$}$$
for all $j\geq2$. In particular, if the maximal subpartition of $\mu$
converges to the trivial partition, then $P_{\mu'}=\{1\}$.
\item The normal closure $\langle(P_\mu)^G\rangle$ of a partition
  subgroup $P_\mu$ of $G$ is the partition subgroup $P_{\mu''}$,
  where $\mu''=(\mu''_2,\mu''_3,\dots)$ is the partition with
$\mu''_j=\max\{i~|~(i,k)\in\mu~,~\forall~k\leq j\}$.
\end{enumerate}
\end{prop}

It is clear that partition subgroups are closed, since any sequence of
elements in a partition subgroup $P_\mu$ which converges in $G$ must
converge in $P_\mu$. The other statements are routine.

\smallskip
Weir obtained specific results for normal partition subgroups of the
finite quotients $G_n(q)$, for all $n,q$, and these also apply to $G$. 

\begin{prop}\cite[Theorem 3]{weir}\label{prop:weir}
Given a normal partition subgroup $P_\mu$, then 
\begin{enumerate}
\item $[P_\mu,G]=P_{\mu'}$, where $|\mu'|$ are all the squares covered
  by $|\mu|$. 
\item Let $P^*$ be the preimage of $Z(G/P_\mu)$ in $G$, then
  $P^*=P_{\hat{\mu}}$ is the normal partition subgroup where
  $|\hat{\mu}|$ are all the squares which do not avoid $\mu$.
\end{enumerate}
\end{prop}

Thus to get the commutator subgroup $[P_\mu,G]$ we ``delete'' the
squares at the corners of $|\mu|$, i.e. if $(i,j)$ and $(i+1,j+1)$ are
both outside $|\mu|$ but $(i,j+1)$ is in $|\mu|$, then we delete
$(i,j+1)$ in $|\mu'|$. On the other hand, if $(i,j)$ and $(i+1,j+1)$ are
both in $|\mu|$ but $(i+1,j)$ is not in $|\mu|$, then we add $(i+1,j)$
to $|\mu|$ to get $|\hat{\mu}|$.

As an example for Proposition~\ref{prop:weir}, we get the
subgroups $\gamma_d(G)$ in the lower central series of $G$ by deleting
successive {\em super diagonals}, where the $d$th super diagonal is the
set of all squares $(i,i+d)\in\N^2$. Thus 
$$\gamma_{d}(G)=P_{(0^{d-1},1,2,3,\dots)}
\qbox{for all $d\in\N$, starting with}\gamma_1(G)=G.$$
Similar considerations allow us to calculate the derived subgroups
$G^{(d)}=[G^{(d-1)},G^{(d-1)}]$ of $G$, starting with 
$$G^{(2)}=[G,G]=\gamma_2(G)=\{1+\sum_{j\geq i+2}a_{ij}e_{ij}\in G\}.$$
Thus elementary commutators calculations (Equation~(\ref{eqn2})), with
$d\geq2$, give 
$$\big[1+\sum_{j\geq i+2^{d-1}}a_{ij}e_{ij},
1+\sum_{j\geq i+2^{d-1}}b_{ij}e_{ij}\big]=
1+\sum_{k\geq i+2^d}c_{ik}e_{ik}.$$
Hence, as partition subgroup,
$$G^{(d)}=P_{(0^{(-1+2^{d-1})},1,2,3,\dots)}\qbox{for all $d\geq2$.}$$
In particular, $G$ is not soluble, because its
derived series does not converge.

Partition subgroups can also be used to show that $G$ is not
hereditarily just infinite.
A profinite group $G$ is {\em hereditarily just infinite} if every every
open subgroup is {\em just infinite} (\cite[Definition
  I.3]{klaas}). That is, every nontrivial closed normal subgroup of any
open subgroup of $G$ has finite
index. By the above discussion, the open subgroups of $G$ have
infinitely many closed normal subgroups of infinite index (e.g. the
subgroups in the lower central series).

\smallskip
From the basic commutator formula~(\ref{eqn2}), we obtain the structure
of the centralisers of partition subgroups of $G$. 

\begin{defn}\label{def:mu-perp}
Let $\mu$ be a partition diagram. The {\em orthogonal partition diagram}
of $\mu$ is the partition diagram
\begin{equation}\label{eqn:mu-perp}
\mu^\perp=\{(k,l)\in\N^2~|~k<l\;\hbox{and}\;
\forall\;(i,j)\in\mu\;\hbox{then}\;k\neq j\;\hbox{and}\;l\neq i\}.
\end{equation}
The {\em centre of $\mu$} is the subpartition diagram
\begin{equation}\label{eqn:centre-mu}
\zeta_\mu=\mu\cap\mu^\perp\qbox{of $\mu$.}
\end{equation}
\end{defn}

For instance, if $\mu=\{(3,4)\}$, then 
$$\mu^\perp=\{(k,l)\in\N^2~|~k\neq4~,~l\neq3\}
\qbox{and}\zeta_\mu=\mu\;,\;\hbox{i.e.,}$$
$$P_{\mu^\perp}=\begin{pmatrix}
1&*&0&*&*&*&*\\
&1&0&*&*&*&*\\
&&1&*&*&*&*\\
&&&1&0&0&0\\
&&&&1&*&\dots *\\
&&&&&\ddots&\dots\end{pmatrix}
\qbox{and}P_{\zeta_\mu}=P_\mu=
\begin{pmatrix}
1&0&0&0&0&0&0\\
&1&0&0&0&0&0\\
&&1&*&0&0&0\\
&&&1&0&0&0\\
&&&&\ddots&\ddots&\dots\end{pmatrix}~.$$

By definition
$\dst C_G(P_\mu)=
\bigcap_{\substack{(i,j)\in\mu\\a_{ij}\in\F_q}}C_G(1+a_{ij}e_{ij})$ and each
$C_G(1+a_{ij}e_{ij})=C_G(1+e_{ij})$ is a partition subgroup of $G$. 
Now $1+ae_{kl}\in C_G(1+e_{ij})$ for each $(i,j)\in\mu$ and $a\in\F_q$
if and only if 
$$1=[1+e_{kl},1+e_{ij}]=1+\delta_{il}e_{kj}-\delta_{jk}e_{il}.$$
So $l\neq i$ and $k\neq j$ for each $(i,j)\in\mu$. In other words,
$$1+e_{kl}\in C_G(P_\mu)
\Longleftrightarrow
(k,l)\in\mu^\perp.$$
Which leads to the following conclusion.

\begin{prop}\label{prop:norm-partition}
Let $\mu$ be a partition diagram.
Then
$$C_G(P_\mu)=P_{\mu^\perp}\qbox{and}
P_{\zeta_\mu}\qbox{is the centre of $P_\mu$.}$$
In particular, $C_G(P_\mu)=\{1\}$ for any open partition subgroup.
\end{prop}


\section{Examples of torsion subgroups}\label{sec:torsion}\
 
The direct limit $\inlim nG_n$ of the $G_n$'s is a discrete torsion
group, and so not a subgroup of $G$. 
Here $\inlim nG_n$ is the group formed by all the square matrices $x$
such that there exists $m\in\N$ for which $x\in G_m$. 
Now, each element of $\inlim nG_n$
can be regarded as a torsion element in $G$ in the obvious way,
by taking each $x\in \inlim nG_n$ to
$\dst\big(xN_n/N_n)_{n\in\N}\in G$. This
mapping, let's call it $\rho$, is an injective homomorphism of abstract
groups, which takes 
$\dst\inlim{n\in\N}G_n$ into $G$ and with the property that its image is
dense in $G$, i.e. $\dst\bigcap_{n\in\N}\im(\rho)N_n/N_n=G$. Note that
$G$ contains ``many'' torsion elements which are not in $\im(\rho)$
(take for instance $1+\dst\sum_{1<j}e_{1j}$).

In~\cite{holub0}, Holubowski studies string subgroups, which form a
large class of torsion discrete subgroups of $G$. In particular string
subgroups cannot contain any open subgroup of $G$.   

\begin{defn}
A matrix $a\in G$ is a {\em string} if $a$ is in the image of some
injective group homomorphism
$$\prod_{n_i>1}G_{n_i}\hookrightarrow G\qbox{of
pairwise diagonal commuting block matrices of size greater than $1$.}$$
Thus $a$ has finite order and $a^{-1}$ is a string with the
same block structure as $a$.
A {\em string subgroup} of $G$ is a subgroup of $G$ formed by strings.
\end{defn}
So a string subgroup $Q$ of $G$ is isomorphic to a subgroup of a
partition subgroup of $G$ of the form $\dst\prod_{n_i>1}G_{n_i}$ for
some non-negative integers $n_i$, for $i\in\N$. 

The equality
$$\prod_{i\in\N}G_{n_i}\cong G/P_\mu\qbox{where}
\mu=(0^{n_1},n_1^{n_2},(n_1+n_2)^{n_3},\dots,
(\sum_{1\leq i\leq j}n_i)^{n_{j+1}},\dots)$$
shows that, regarded as abstract groups (as in \cite{holub0}), string
subgroups are the complements of normal partition subgroups. That is,
$$G=P_\mu\cdot\prod_{i\in\N}G_{n_i}\qbox,
P_\mu\cap\bigg(\prod_{i\in\N}G_{n_i}\bigg)=\{1\}
\qbox{and}P_\mu\trianglelefteq G.$$


\section{Free subgroups of $G$}\label{sec:free}

In this section, we let $q=p$, and $G=\prlim{n\in\N}G_n(p)$.

We discuss a particular discrete subgroup of $G$ investigated by
Holubowski (cf.~\cite{holub0}), and we also look at its closure in $G$. 
This subgroup of $G$ is the product of two string
subgroups, but is not a string subgroup itself. 

\begin{defn}\label{def:free}
Let
\begin{align*}
s&=1+\sum_{n\in\N}e_{2n-1,2n}=\begin{pmatrix}
1&1\\&1\\&&1&1\\&&&1\\&&&&\ddots\end{pmatrix}
\qbox{and}\\
t&=1+\sum_{n\in\N}e_{2n,2n+1}=\begin{pmatrix}
1\\&1&1\\&&1&\\&&&1&1\\&&&&1\\&&&&&\ddots\end{pmatrix}
\end{align*}
and regard $s$ and $t$ as elements of $G$.

Let $F=\langle x\rangle*\langle y\rangle$ be the free
product of two groups of order $p$. 
\end{defn}

Holubowski (\cite[Theorem~1]{holub0}) defines a function 
$\varphi~:~F\to G$, by  
$$\varphi(x)=s=1+\sum_{n\in\N}e_{2n-1,2n}\qbox{and}
\varphi(y)=t=1+\sum_{n\in\N}e_{2n,2n+1}$$
and proves that $\varphi$ is an injective group homomorphism. 
\cite[Theorem~1]{holub0} shows that the image is contained in the
intersection of the subgroup of so-called {\em banded matrices} with the
subgroup of periodic matrices of $G$. A {\em banded matrix} is a matrix
$(a_{ij})\in G$ for which there exists $d\in\N$ such that $a_{ij}=0$
whenever $j>i+d$. In particular, $\im(\varphi)$ is not a closed subgroup
of $G$ because $\im(\varphi)<\bigcap_n\im(\varphi)N_n/N_n$. 
We let $Q=\ovl{\varphi(F)}$ be the closure of
$\varphi(F)$ in $G$.

A word $w(x,y)=x^{a_1}y^{b_1}\cdots x^{a_l}y^{b_l}$ is mapped to 
\begin{align*}
\varphi(w(x,y))=w(s,t)&=s^{a_1}t^{b_1}\cdots s^{a_l}t^{b_l}\\
&=1+\sum_{n\in\N}\left(\sum_{0\leq j<2l}A_je_{2n-1,2n+j}+
\sum_{1\leq j<2l}B_je_{2n,2n+j}\right)
\end{align*}
where the coefficients $A_j,B_j$ are monomials in
$a_1,\dots,a_l,b_1,\dots,b_l$ of degrees at most $j+1$ and $j$
respectively. 
For instance, for any $a,b,c,d\in\F_p$, 
\begin{equation}\label{eqn:word}
\begin{array}{l}
s^at^bs^ct^d=\\
1+\dst\sum_{n\in\N}\bigg(
(a+c)e_{2n-1,2n}+(ab+cd+ad)e_{2n-1,2n+1}+\\
\hspace{1cm}+abce_{2n-1,2n+2}+abcde_{2n-1,2n+3}+\\
\hspace{1cm}+(b+d)e_{2n,2n+1}+bce_{2n,2n+2}+bcde_{2n,2n+3}\bigg)\\
=\begin{pmatrix}
1&a+c&ab+cd+ad&abc&abcd&0&\dots\\
0&1&b+d&bc&bcd&0&\dots\\
0&0&1&a+c&ab+cd+ad&abc&\dots\\
&&&\ddots&\ddots\end{pmatrix}~.
\end{array}
\end{equation}

Recall from Definition~\ref{def:periodic} that $X[2]$ is the matrix
obtained from $X\in G$ by deleting the first $2$ rows and columns.
From Equation~(\ref{eqn:word}) and elaborating by induction on it, we
record the following. 

\begin{prop}\label{prop:f-periodic}
For any $w(x,y)\in F$, we have
$$\varphi(w(x,y))=w(s,t)=w(s,t)[2].$$
Moreover, the length of $w(x,y)$ can be read from the last nonzero
squares in the first two rows of
$w(s,t)=x^{a_1}y^{b_1}\cdots x^{a_l}y^{b_l}$. Namely, 
\begin{itemize}
\item[(i)] if $a_1b_l\neq0$, i.e. $a_j,b_j\neq0$ 
for all $1\leq j\leq l$, then the last nonzero squares in the first two
rows of $w(s,t)$ are $(1,2l+1)$ and $(2,2l+1)$ respectively; 
\item[(ii)] if $a_1=0\neq b_l$, (i.e. $a_1$ is the only zero exponent) 
then the last nonzero squares in the
  first two rows of $w(s,t)$ are $(1,2l-1)$ and $(2,2l+1)$;
\item[(iii)] if $a_1=b_l=0$ and no other exponent is zero, then the last
  nonzero squares in the first two rows of $w(s,t)$ are $(1,2(l-1))$ and
  $(2,2l)$. 
\end{itemize}
In particular, we obtain $2$-periodic elements of $G$ whose last
nonzero squares in any two successive rows $(i,j),(i+1,k)$ are such that
$k-j\in\{0,2\}$. 
Therefore 
$$Q<\{X\in G~|~X=X[2]\}$$
is a closed subgroup of infinite index in the subgroup of $2$-periodic
elements of $G$. Furthermore,
$$\varphi^{-1}(\gamma_d(G))=\gamma_d(F)\qbox{for all $d\in\N$.}$$
\end{prop}

A counting exercise gives the indices $|PN_n:N_n|=p^{2n-3}$ and 
$|QN_n:N_n|=p^{\xi_n}$, where $P$ is the subgroup of $2$-periodic
elements of $G$ and $\xi_n=\dst n-2+\lfloor\frac{n+1}2\rfloor$, this
latter obtained inductively on $n$.

\begin{remark}\

\begin{enumerate}
\item It is important to emphasise that Holubowski regards $\varphi(F)$
  (as most of the subgroups he investigates in
  \cite{holub1,holub0,holub2}) as a discrete group, and this can be seen
  from the fact that $\im(\varphi)$ is not closed in $G$. 
Recall that $F$ is a pro-$p$ group for the topology defined by taking 
the set $\mathcal F$ of all the subgroups of $F$ of finite $p$-power index
(cf.~\cite[Example (iii), p. 29]{hig}). Letting $R$ run
through all the normal subgroups of $F$ of finite $p$-power index, we
obtain all the $2$-generated finite $p$-groups as the quotients $F/R$. 
For instance $\gamma_2(F)\in\mathcal F$ and
$F/\gamma_2(F)\cong C_p\times C_p$. 

\item Let us also mention the description of $F$ from
\cite[p. 28]{serre}, i.e. that of the fundamental group of a tree with
fundamental domain a segment 
$\xymatrix{\langle x\rangle\ar@{-}[rr]^1&&\langle y\rangle}$.
\end{enumerate}
\end{remark}


\section{Nottingham group}\label{sec:nott}
As pointed out in the concluding remark of \cite{holub1}, the Nottingham
group can be seen as a subgroup of $G$. We make this embedding 
of topological groups explicit in this section. There are equivalent
definitions of the Nottingham group. We follow \cite{c2}.

\begin{defn}
The Nottingham group $\CN=\CN_q$ is the group of algebra automorphisms
of $\F_q[[t]]$ of the form
$$t\mapsto t+\sum_{j\geq2}a_jt^j\quad a_j\in\F_q.$$
\end{defn}
R. Camina investigated the subgroups of $\CN$ and proved that $\CN$
contains every countably based pro-$p$ group as a closed subgroup.

Now the tantalising fact that $G=\prlim nG_n$ is countably based as
pro-$p$ group, implies by Camina's result that $G$ embeds into
$\CN$ as a closed subgroup. On the other hand, by linear
algebra, the elements of $\CN$ can be expressed as
infinite unitriangular matrices, i.e. elements of $G$, and therefore
$\CN$ is a subset of $G$; but it certainly cannot be the whole of $G$,
because $\CN$ only consists of algebra automorphisms.

The convention is that $\CN$ acts on the right of $\F_q[[t]]$. So, we
can identify the nonconstant elements of $\F_q[[t]]$ as infinite row vectors
$$\sum_{j\geq1}a_jt^j\in\F_q[[t]]\longleftrightarrow
(a_1,a_2,a_3,\dots)\in\F_q^{\aleph_0}.$$
Matrix multiplication induces a linear transformation of $\F_q^{\aleph_0}$,
$$(a_1,a_2,a_3,\dots)\mapsto(a_1,a_2,a_3,\dots)x=
(a_1,a_2+x_{12}a_1,\dots,a_j+\sum_{1\leq i<j}x_{ij}a_{i},\dots)~,$$
for all $(x_{ij})\in G$.
which translates as function on $\F_q[[t]]$ as follows:
$$\sum_{j\geq1}a_jt^j\mapsto 
a_1t+(a_2+x_{12}a_1)t^2+\dots+\left(a_j+\sum_{1\leq i<j}x_{ij}a_{i}\right)t^j+\dots$$

Given that 
$\CN=\langle e_1[\alpha_c],e_2[\alpha_c]~|~1\leq c\leq f\rangle$, 
where 
$\alpha_1,\dots,\alpha_f\in\F_q$ generate $\F_q$ as $\F_p$-vector
space, and 
$$e_r[\alpha_c]~:~t\mapsto t+\alpha_ct^{r+1}\in\CN~,\qbox{we have}$$
$$(\sum_ja_jt^j)e_r[\alpha_c]=\sum_ja_j(t+\alpha_ct^{r+1})^j=
\sum_ja_j\bigg(\sum_{0\leq i\leq j}\binom ji\alpha_c^it^{ri+j}\bigg).$$
This suggests the following mapping $\sigma:\CN\to G$, defined on the
generators of $\CN$ by
\begin{equation}\label{eqn:sigma}
\sigma(e_r[\alpha_c])=g_r[\alpha_c]\stackrel{\hbox{\scriptsize{def}}}=
\sum_{1\leq i\leq j}
\binom i{\frac{j-i}r}\alpha_c^{\frac{j-i}r}~e_{ij}=
1+\sum_{1\leq i<j}
\binom i{\frac{j-i}r}\alpha_c^{\frac{j-i}r}~e_{ij}~,
\end{equation}
where the sums run over all the positive integers $i\leq j$,
resp. $i<j$, such that 
$$\frac{j-i}r\in\Z\qbox{and}j\leq 2i$$
and all the other coefficients are zero. 
Note that in the first sum
$g_r[\alpha]_{ii}=1$ for all $i\in\N$. 

In particular, for $\alpha_1=1$, if we write $e_r=e_r[1]$ and
$g_r=g_r[1]$, then the $i$th row of $g_r$ contains the $i$th row of
Pascal's triangle starting from the diagonal $1$, and spaced by $(r-1)$
zeroes between each coefficient in a row. 
Note that $g_r[\alpha_c]\in\gamma_r(G)$ for all $r\geq1$.

For example, 
$$g_1=\begin{pmatrix}1&1&0&\dots\\
&1&2&1&0&\dots\\
&&1&3&3&1&0&\dots\\
&&&\ddots&\ddots&\ddots\end{pmatrix}
\qbox{and}$$
$$g_2=\begin{pmatrix}
1&0&1&0&\dots\\
&1&0&2&0&1&0&\dots\\
&&1&0&3&0&3&0&1&0\dots\\
&&&&&&\ddots&\ddots&\ddots\end{pmatrix}$$
Matrix multiplication yields 
$$(a_1,a_2,\dots)g_1=(a_1,a_1+a_2,2a_2+a_3,a_2+3a_3+a_4,\dots)$$
which corresponds to
$$a_1t+(a_1+a_2)t^2+(2a_2+a_3)t^3+(a_2+3a_3+a_4)t^4+\dots\in\F_q[[t]]$$
and so gives in particular $(1,0,0,\dots)g_1=(1,1,0,\dots)$,
i.e. $te_1=t+t^2$. Accordingly, for any ``canonical'' vector
$t^i\in\F_q[[t]]$ the corresponding ``canonical'' row vector
$v_i\in\F_q^{\aleph_0}$ has a unique nonzero coefficient equal to $1$ in
the $i$th coordinate, so that $v_ig_r[\alpha_c]$ is the $i$th row of
$g_r[\alpha_c]$,
$$(0^{i-1},1,0^{r-1},i\alpha_c,0^{r-1},\binom i2\alpha_c^2,
0^{r-1},\dots,0^{r-1},\alpha_c^i,0^{r-1},\dots)$$
which corresponds to the element
$$t^i+i\alpha_ct^{i+r}+\alpha_c^2\binom i2t^{i+2r}+
\dots+\alpha_c^it^{i(r+1)}\in\F_q[[t]].$$

Routine computations give
$$(g_1^2)_{ij}=\sum_{k\geq1}(g_1)_{ik}(g_1)_{kj}=
\sum_{i\leq k\leq2i}\binom i{k-i}\binom k{j-k}$$
for\; $i\leq j\leq4i$\;
and $(g_1^2)_{ij}=0$ otherwise. That is, a ``row-palindrome'' matrix
$$\begin{pmatrix}
1&2&2&1&\dots\\&1&4&8&10&8&4&1\dots\\&&\dots\end{pmatrix}.$$
Since each row of $g_r[\alpha_c]$ has finitely many nonzero entries, a
recursive algorithm (or a more elaborate procedure) gives us the
inverse; for instance 
\begin{align*}
(g_1)^{-1}&=\begin{pmatrix}
1&-1&2&-5&14&\dots\\
&1&-2&5&14&\dots\\
&&1&-3&9&\dots\\
&&&1&-4&\dots\\
&&\ddots&&\ddots\end{pmatrix}\qbox{and}\\
(g_2)^{-1}&=\begin{pmatrix}
1&0&-1&0&3&0&-12\dots\\
&1&0&-2&0&7&0\dots\\
&&1&0&-3&0&12\dots\\
&&&1&0&-4&0\dots\\
&&\ddots&&\ddots\end{pmatrix}~.
\end{align*}
In particular, each row of $g_r[\alpha_c]^{-1}$ has infinitely many
nonzero terms, and $g_r[\alpha_c]^{-1}\in\gamma_r(G)$ for all $r\geq1$.

The key point is that the elements in $\im(\sigma)$ are entirely
determined by their first row, where $\sigma$ is defined by
Equation~(\ref{eqn:sigma}).
That is, if $x\in\CN$ is given by
$tx=t+\sum_ja_jt^j$, then the equation $(t^i)x=(t+\sum_ja_jt^j)^i$
defines the coefficients in the $i$th row of $\sigma(x)$.

It is routine to check that the matrices $g_r[\alpha_c]$ correspond to
the image under $\sigma$ of the algebra automorphisms
$e_r[\alpha_c]\in\CN$, and that they are subject to the same relations.

\begin{prop}\label{prop:cn}
$\sigma(\CN)$ is a closed subgroup of $G$ of infinite index in $G$.
In particular, $\sigma(\CN)$ does not contain any open subset of $G$.
\end{prop}

\begin{proof}
We have seen above that $\sigma$ is a homomorphism of abstract groups,
and it is clearly injective. 
For any $g\in\sigma(\CN)$ and for any $\delta\in(0,1)$, the
open ball $B(g,\delta)$ is not contained in $\sigma(\CN)$, since it
contains infinitely many linear transformations of $\F_q^{\aleph_0}$ which are
not algebra automorphisms.
Therefore, $\sigma(\CN)$ does not contain any open subset of $G$ and has
infinite index.
To prove that $\sigma$ is continuous, and so that $\sigma(\CN)$ is a
closed subgroup of $G$, we show that the preimage by
$\sigma$ of any neighbourhood 
$B(\sigma(u),\epsilon^{n-1})=\sigma(u)N_n$
is a neighbourhood of $u$ for any $u\in\CN$ and $n\in \N$.
Note that $\sigma^{-1}(N_n)=
\langle e_m[\alpha_c]~|~1\leq c\leq f~,~m\geq n\rangle$ is an open
normal subgroup of $\CN$ (\cite{c2}). So we have
$$\sigma^{-1}\big(\sigma(u)N_n\big)=
u~\langle e_m[\alpha_c]~|~1\leq c\leq f~,~m\geq n\rangle$$
which is an open set of $\CN$, as required.
\end{proof}

\smallskip
Next, we turn to the normaliser $N_G(\sigma(\CN))$ of 
$\sigma(\CN)$ in $G$. 
Klopsch proved in \cite{klopsch} that every automorphism of $\CN$ is
standard, provided $p\geq5$. That is, 
$\Aut(\CN)\cong\mathcal A(q)\rtimes\Aut(\F_q)$, where $\mathcal A(q)$ is
the group of all the algebra automorphisms
$\{t\mapsto\sum_{n\geq1}\lambda_nt^n~|~\lambda_n\in\F_q~,~\lambda_1\neq0\}$
of $\F_q[[t]]$. 

\begin{prop}\label{nott:norm}
Suppose $p\geq5$. Then $N_G(\sigma(\CN))=\sigma(\CN)$,
i.e. $\sigma(\CN)$ is selfnormalising in $G$. 
\end{prop}

\begin{proof}
Consider the inclusion $N_G(\sigma(\CN))/C_G(\sigma(\CN))\hookrightarrow
\Aut(\sigma(\CN))\cong\Aut(\CN)$ given by
mapping $gC_G(\sigma(\CN))$ to conjugation by $g$ in $\sigma(\CN)$. From
the elementary commutator relation~(\ref{eqn2}), we observe that
$C_G(\sigma(\CN))=\cap_{r,c}C_G(g_r[\alpha_c])=\{1\}$, where $r\in\N$ and 
$1\leq c<f$. 
Moreover, $G$ is a pro-$p$ group and so 
$N_G(\sigma(\CN))\hookrightarrow S_p$, where $S_p$ is a Sylow $p$-subgroup of
$\Aut(\CN)\cong\mathcal A(q)\rtimes\Aut(\F_q)$. Now, 
$\Aut(\F_q)=\langle\Phi~:~\alpha\mapsto\alpha^p\rangle$ is the cyclic
group spanned by the Frobenius homomorphism $\Phi$, which has order
$f=\log_pq$ (\cite[VII.5~Theorem~12]{lang}). In particular, for
$\alpha\in\F_q-\F_p$, we have 
$\Phi(g_1[\alpha])\not\equiv g_1[\Phi(\alpha)]\pmod{\gamma_2(G)}$,
implying that this mapping cannot be given by conjugation by an element
of $G$. It follows that $N_G(\sigma(\CN))$ is isomorphic to a
$p$-subgroup of $\mathcal A(q)$.
Since $\CN$ is the unique Sylow $p$-subgroup of 
$\mathcal A(q)\cong\CN\rtimes\F_q^\times$, the result follows.
\end{proof}


\section{Field extensions}\label{sec:field-ext}

Given $q=p^f$ for $f\in\N$ and $p$ an odd prime, let us regard $\F_q$ as
an $f$-dimensional $\F_p$-vector space.
Left multiplication in $\F_q$ induces an injective ring homomorphism
between endomorphism rings of vector spaces 
$$\wh{\alpha}_f~:~\End\F_q\to\End(\F_p^f)\cong\Mat_f(\F_p).$$
Choosing such $\wh{\alpha}_f$ induces an injective group homomorphism
(cf. Definition~\ref{def:g})
$$\alpha_f~:~G(q)\to G(p).$$
Note that $\wh\alpha_f(\lambda)\neq0\in\Mat_f(\F_p)$ if and only if
$\lambda\neq0\in\F_q$, in which case $\wh\alpha_f(\lambda)\in\GL_f(p)$. 
Also, $\wh\alpha_f(\lambda)=\lambda I_f$ if and only if $\lambda$ is in
the subfield $\F_p$ of $\F_q$.

\begin{lemma}\label{lem:alpha-cont}
The map $\alpha_f$ is continuous. So $G(q)$ is isomorphic to a
closed subgroup of $G(p)$.
\end{lemma}

\begin{proof}
Write $\alpha=\alpha_f$, and $N_k(p)$ and $N_k(q)$ for the
normal open subgroups of $G(p)$ and $G(q)$ respectively.
Let $\alpha(x)\in G(p)$ and $n\in\N$. 
Then
\begin{align*}
\alpha^{-1}\big(B(\alpha(x),\epsilon^n)\big)&=
\alpha^{-1}\big(
\{z\in G(p)~|~z^{-1}\alpha(x)\in N_{n+1}(p)\}\big)\\
&\stackrel{(\dagger)}=\alpha^{-1}\big(\{\alpha(y)~|~y\in
G(q)~,~\alpha(y)^{-1}\alpha(x)\in N_{n+1}(p)\}\big)\\
&=\{y\in G(q)~|~\alpha(y^{-1}x)\in N_{n+1}(p)\}\\
&=\{y\in G(q)~|~y^{-1}x\in N_k(q)\}\\
&=B(x,\epsilon^k)
\end{align*}
where $k\in\N$ is defined by the inequalities
$kf\leq n+1<(k+1)f$ and $(\dagger)$ follows from taking only those 
$z\in G(p)$ which have nonempty preimage by $\alpha$. 
Therefore, the preimage under $\alpha$ of an open
neighbourhood of $\alpha(x)$ is an open neighbourhood of $x$ for all
$x\in G(q)$, which proves that $\alpha$ is continuous
(\cite[Ch.~3,~Theorem~1]{kelley}).  
\end{proof}

In particular, for any $f,k\in\N$, we have
$\alpha_f^{-1}(N_{kf}(p))=N_k(q)$.

\smallskip
Conversely, the ring inclusion $\F_p\hookrightarrow\F_q$ induces a
continuous injective group homomorphism 
$$\beta_f~:~G(p)\to G(q)\qbox{and we have}$$
$$\alpha_f\beta_f(x)=x\otimes I_f=
\left(\begin{array}{r|r|r|r}
I_f&x_{12}I_f&x_{13}I_f&\dots\\
\hline
&I_f&x_{23}I_f&\dots\\
\hline
&&I_f&\dots\\
\hline
&&&\ddots\end{array}\right)\in G(p)\qbox{for all $x\in G(p)$.}$$

For short, fix $f\in\N$ and let $\alpha=\alpha_f$,
$\wh\alpha=\wh\alpha_f$ and $\beta=\beta_f$. 
Let $H=\im(\alpha)\leq_cG(p)$ and $K=\im(\beta)\leq_cG(q)$. 
Clearly neither $H$ nor $K$ are open subgroups because they have
infinite index in $G(p)$ and $G(q)$ respectively.
To find the index of $H$ in $G(p)$, we regard the elements of $G(p)$ in
$f\times f$ block form, where each block is of the form
\begin{equation}\label{eqn:f-block}
\xymatrix{\big((f-1)i+1,(f-1)j+1\big)\ar@{.}[rr]\ar@{.}[dd]
\ar@{.}[rrdd]&&
\big((f-1)i+1,fj\big)\ar@{.}[dd]\ar@{.}[lldd]\\
{}\ar@{.}[rr]&&{}\\
\big(fi,(f-1)j+1\big)\ar@{.}[rr]&&\big(fi,fj\big)}
\qbox{with $i\leq j$.}
\end{equation}
In each of these blocks with $i<j$, the image of $\wh\alpha$
is isomorphic to a copy of $\F_q$, so that in each block, the index is
equal to $|\Mat_f(\F_p)~:~\wh\alpha(\F_q)|=p^{f^2-f}$. There are
countably infinitely many such blocks, giving 
$|G(p):H|=(p^{f^2-f})^{\aleph_0}$.

Similarly, to calculate the index of $G(p)$ in $G(q)$ we see that for
each coefficient $(i,j)$ with $i<j$, we have an index
$|\F_q:\F_p|=p^{f-1}$, so that $|G(q):K|=(p^{f-1})^{\aleph_0}$. 

\smallskip
From Proposition~\ref{prop:bier-weir}, we gather that 
$H\not\trianglelefteq G(p)$ and that $K\not\trianglelefteq G(q)$ since
they are not partition subgroups.

In order to determine the normalisers $N_{G(p)}(H)$ and $N_{G(q)}(K)$,
we use a result of Weir on the automorphisms of the finite groups
$G_n(q)$. For convenience, we have put Weir's theorem and some technical
considerations of conjugation in Appendix~\ref{app:aut-gn}.

\begin{prop}\label{prop:field-ext-norm}
Suppose the above notation. The following hold.
\begin{itemize}
\item[(i)] $N_{G(q)}(K)=K$, where $K=\beta_f(G(p))$.
\item[(ii)] $N_{G(p)}(H)=H$, where $H=\alpha_f(G(q))$.
\end{itemize}
\end{prop}

\begin{proof}
We prove the second part: $N_{G(p)}(H)=H$. 
Let $\alpha=\alpha_f$ and $H_m=HN_m(p)/N_m(p)$ for any $m\in\N$.
We first show that in the finite quotients $G_{nf}(p)$ we have
$$N_{G_{nf}(p)}(H_{nf})=H_{nf}C_{G_{nf}(p)}(H_{nf})$$
where 
$C_{G_{nf}(p)}(H_{nf})=
\langle1+e_{ij}~|~1\leq i\leq f~,~n-f<j\leq n\rangle$ 
is the subgroup of all elements of
$G_{nf}(p)$ whose only nonzero nondiagonal squares are in the upper
right $f\times f$ corner, by Lemma~\ref{lem:centr-h} below.  

Using Theorem~\ref{thm:weir-aut}, it remains to show that no other
automorphism of $H_{nf}$ of $p$-power order can be expressed as a
conjugation by an element of $G_{nf}(p)$. 
Thus we need to consider $\mathcal P$ and possibly
$\mathcal L$ in case $p$ divides $f=\log_pq$. 
Now, a field automorphism of $G_n(q)$ becomes an automorphism of
$H_{nf}$ which fixes all the elements
$\alpha(1+e_{r,r+1})$ and therefore cannot be given by a conjugation by
a matrix in $G_{nf}(p)$ because the elements which centralise
$\alpha(K)N_{nf}(p)/N_{nf}(p)$ also centralise $H_{nf}$. So
Lemma~\ref{lem:centr-h} proves that a field automorphism on $H_{nf}$
cannot be given by an inner automorphism of $G_{nf}(p)$. 

Finally, neither central, nor extremal automorphisms of $H_{nf}$ can be
given by inner automorphisms of $G_{nf}(p)$ by a similar argument to
that used in the proof of \cite[Theorem 8]{weir}. Indeed, in
Equation~(\ref{eq:conjug}) in the proof of Lemma~\ref{lem:centr-h}
below, if a morphism 
$\alpha_f(1+e_{r,r+1})\mapsto\alpha_f(1+e_{r,r+1}+ce_{1,n})$, for $r>1$
and $c\in\F_q$, were given by an inner automorphism of $G_{nf}(p)$, say
conjugation by $1+\sum_{i<j}a_{ij}e_{ij}$, then 
Equation~(\ref{eq:conjug}) has no solution; that is, on the one hand we
would need $a_{ij}=0$ for all $f<i<j\leq n$, and on the other,
conjugating $\alpha_f(1+e_{r,r+1})$ by such element is not of the form 
$\alpha_f(1+e_{r,r+1}+ce_{1,n})$.

It follows that $N_{G_{nf}(p)}(H_{nf})=H_{nf}C_{G_{nf}(p)}(H_{nf})$ as
asserted. Now, to obtain the normaliser $N_{G(p)}(H)$, we let
$n\to\infty$, and the claim follows.

The proof that $N_{G(q)}(K)=K$ can be handled in a similar way and
follows easily from the above. We leave this to the reader.
\end{proof}


\section{$p$-adic variation}\label{sec:p-adic}

In this section, we consider a variant on the pro-$p$ groups
$G(q)=\prlim{n\in\N}G_n(q)$ of Definition~\ref{def:g}. 
Namely, let $p$ be an odd prime and for any $n,k\in\N$, write
$$G_n(\Z/p^k)=\{x\in\GL_n(\Z/p^k)~|~x_{ij}=0\;\forall\;i>j\qbox{and}
x_{ii}=1\}\qbox{for all $(n,k)\in I$,}$$
where $I=\{(n,k)\in\N^2\}$ is a poset for the order relation 
$(n,k)\leq (m,l)$ if and only if $n\leq m$ and $k\leq l$. Thus $I$ is a
directed system.

For any $(n,k)\leq(m,l)$ in $I$, there are obvious surjections:
\begin{itemize}
\item $G_m(\Z/p^l)\to G_{n}(\Z/p^l)$ analogous to the surjective group
  homomorphisms from Section~\ref{sec:intro}, using that 
$G_{n+1}(\Z/p^l)\cong(\Z/p^l)^n\rtimes\GL_n(\Z/p^l)$, where
  $(\Z/p^l)^n$ is the natural module for $\GL_n(\Z/p^l)$, and 
\item $G_m(\Z/p^l)\to G_m(\Z/p^k)$ induced by the reduction
  $\Z/p^l\to\Z/p^k$ of the coefficients. 
\end{itemize}
Hence we get an inverse system in which all the maps are surjective
$$\xymatrix{
&\vdots\ar[d]&\vdots\ar[d]&\vdots\ar[d]\\
\dots\ar[r]&G_{n+1}(\Z/p^{k+1})\ar[d]\ar[r]&G_{n}(\Z/p^{k+1})\ar[d]\ar[r]&
G_{n-1}(\Z/p^{k+1})\ar[d]\ar[r]&\dots\\
\dots\ar[r]&G_{n+1}(\Z/p^k)\ar[d]\ar[r]&G_{n}(\Z/p^k)\ar[d]\ar[r]&
G_{n-1}(\Z/p^k)\ar[d]\ar[r]&\dots\\
\dots\ar[r]&G_{n+1}(\Z/p^{k-1})\ar[d]\ar[r]&G_{n}(\Z/p^{k-1})\ar[d]\ar[r]&
G_{n-1}(\Z/p^{k-1})\ar[d]\ar[r]&\dots\\
&\vdots&\vdots&\vdots}$$
The inverse limit of this system is the group
$$G(\Z_p)=\prlim{(n,k)\in I}G_n(\Z/p^k)
\cong\prlim nG_n(\Z_p)\cong\prlim kG(\Z/p^k),$$
where $\Z_p$ denotes the ring of $p$-adic integers.
By definition, $G(\Z_p)$ is a pro-$p$-group
(\cite[Proposition~2.2.1]{RZ} or \cite[Theorem~1.2.5 (b)]{wilson}),
because the class of pro-$p$ groups is closed under taking closed
subgroups and arbitrary direct products. It follows from the product and
subgroup topologies that the open sets of $G$ must be the cosets of
the factor groups $G/U_n(k)$ where  
\begin{align*}
U_{n}(k)&=\{x\in G(\Z_p)~|~x_{ij}\in p^k\Z_p\;\forall\;1\leq i<j\leq n\}\\
&=\ker\big(\xymatrix{G(\Z_p)\ar@{->>}[r]&G_n(\Z/p^k)}\big)\\
&=\left\{\left(\begin{array}{r|r}G_n(p^k\Z_p)&*\\
\hline0&*\end{array}\right)~\bigg|~*\in\Z_p\right\}
\end{align*}
for all $n,k\in\N$.
So $U_n(k)\triangleleft G$ and $G/U_n(k)\cong G_n(\Z/p^k)$.

For later use, we want to extract from this set a
{\em filtration}, i.e. a totally ordered set of open normal subgroups of
$G$. 
For each $n\in\N$, let 
\begin{equation}\label{eqn:p-adic-base}
V_n=\left\{\left(\begin{array}{r|r}\dst
  G_n(p^n\Z_p)&\hbox{\Large*}\\
\hline
\hbox{\Large0}&\hbox{\Large*}\end{array}\right)~\bigg|~*\in\Z_p\right\}
\end{equation}
so that $V_n\trianglelefteq_o G$, with $V_1\geq V_2\geq V_3\geq\dots$
and $\dst G/V_n\cong G_n(\Z/p^n)$ for all $n\in\N$. 

To prove that this is a base of open neighbourhood of $1$ in
$G$, it suffices to show that each $U_n(k)$ is a union of cosets of the
$V_m$'s. By inspection, we obtain
$$U_n(k)=\bigcup_{g\in G_m(p^k\Z_p/p^m\Z_p)}
\left(\begin{array}{r|r}
g&\hbox{\Large*}\\\hline
\hbox{\Large0}&\hbox{\Large*}\end{array}\right)
\qbox{where}
m=\max\{n,k\},\qbox{and}$$
the coefficients of $g\in G_m(p^k\Z_p/p^m\Z_p)$ are the coefficients of
$g$ in $p^m\Z_p$ running through a set of
representatives of $p^k\Z_p/p^m\Z_p\cong\Z/p^{(m-k)}$. 
That is,
$$U_n(n)=V_n\qbox,
U_n(k)=\bigcup_{g\in G_n(p^k\Z_p/p^n\Z_p)}\left(\begin{array}{r|r}
g&\hbox{\Large0}\\\hline
\hbox{\Large0}&\hbox{\Large I}_\infty\end{array}\right)V_n
\qbox{if $k<n$, and}$$
$\dst U_n(k)=\bigcup_{g\in G_k(p^n\Z_p/p^k\Z_p)}\left(\begin{array}{r|r}
g&\hbox{\Large0}\\\hline
\hbox{\Large0}&\hbox{\Large I}_\infty\end{array}\right)V_k$\quad if
$k>n$.

\smallskip
Although $G(\Z_p)$ is defined over the $p$-adic integers, it is clear
from the definition of a $p$-adic analytic group (\cite[Section~9]{dix})
that $G(\Z_p)$ is not $p$-adic analytic because one cannot find 
homeomorphisms between the open subsets of $G(\Z_p)$ and
$\Z_p$-modules of finite rank.
However, for each $n\in\N$, the group $G_n(\Z_p)$ is the
prototype of a compact $p$-adic analytic group (cf. \cite[\S~5.1]{dix}).
From the fact that $G(\Z_p)\cong\prlim nG_n(\Z_p)$, we observe that the
inverse limit of compact $p$-adic analytic groups need not be $p$-adic
analytic. 

\smallskip
Obviously, there are some similarities between subgroup structures of
$G(\Z_p)$ and $G(\F_q)$ for any $q$. In particular, we have partition
subgroups and free products (of the form $\Z_p*\Z_p$). 
Using the ideals $p^k\Z_p$ in $\Z_p$ we note that $G(\Z_p)$ has in fact
a plethora of closed subgroups. We flag some ``ideal'' partition
subgroups. 

\begin{prop}\label{prop:p-adic-part}
Let $I=p^k\Z_p$ be an ideal of $\Z_p$ and $\mu$ a partition diagram. Let
$$P_\mu(I)=\langle1+p^ke_{ij}~|~(i,j)\in\mu\rangle\leq G(\Z_p).$$
The following hold.
\begin{enumerate}
\item[(i)] $P_\mu(I)$ is closed.
\item[(ii)] $P_\mu(I)$ is open if and only if $I=\Z_p$ and $\mu$
  converges to a partition such that there exists $N\in\N$ with
  $(i,j)\in\mu$ of all $1\leq i<j$ and all $j\geq N$.
\item[(iii)] $P_\mu(I)$ is normal in $G(\Z_p)$ if and only if $P_\mu(I)$
  is characteristic in $G(\Z_p)$ if and only if the partition subgroup
  $P_\mu$ is normal in $G(\Z_p)$.
\end{enumerate}
\end{prop}

The proof is straightforward using Proposition~\ref{prop:bier-weir}. 


\section{Hausdorff dimension of closed subgroups of $G$}\label{sec:hausdorff}

In this section, we apply \cite[Proposition~2.6]{aber} to some closed
subgroups of $G$ and calculate their Hausdorff dimension. Note that in
\cite{aber}, the author refers to the {\em Billingsley dimension}
instead of {\em Hausdorff}, which is defined over the set of real
numbers. We adopt the terminology used in later papers (\cite{bk,bs}),
which use Abercrombie's results too. We refer
the reader to \cite{falconer} for an in-depth background on fractal
dimensions and measure theory.
We limit ourselves to the essential facts as they apply to $G=G(q)$
from Definition~\ref{def:g}, and include an appendix with some
additional theory which may be useful to the reader. 
For convenience, we take as definition of Hausdorff dimension that given
in Abercrombie's result.

\begin{defn}\label{def:hausdorff}\cite[Proposition~2.6]{aber}
Let $H=\dst\prlim{n\in\N}H_n$, be a closed subgroup of 
$G\dst=\prlim n~G_n=\prlim n ~G/N_n$, where $H_n=HN_n/N_n$, and 
where $N_1\geq N_2\geq N_3\geq\dots$ is a filtration of $G$ by open
normal subgroups. 
The {\em Hausdorff dimension of $H$} is the real number 
\begin{equation}\label{prop:aber}
\begin{array}{l}
\dim(H)=\dst\lim_{n\to\infty}\frac{\log|H_n|}{\log|G_n|}
\qbox{whenever the limits exists.}\\
\hbox{Otherwise}\quad
\dim(H)\geq\dst\liminf{n\to\infty}\frac{\log|H_n|}{\log|G_n|}
\end{array}
\end{equation}

The {\em Hausdorff spectrum} of $G$ is the subset 
$$\spec(G)=\{\dim(H)~|~H\leq_cG\}$$
of the Hausdorff dimensions of all the closed subgroups of $G$. Thus 
$\{0,1\}\subseteq \spec(G)\subseteq[0,1]$.
\end{defn}

Recall from the elementary law of logarithms 
$$e^{\ln a}=a=b^{\log_ba}=e^{\ln b\log_ba}\qbox{that}
\log_ba=\dst\frac{\ln a}{\ln b},$$
so that we can take any
base for the logarithm defining $\dim(H)$. We will take $\log=\log_q$
unless otherwise specified.

From the definition, if $|H|<\infty$, then $\dim(H)=0$, and consequently, any
subset $K\subseteq G$ of finite index has $\dim(K)=1$. Therefore, the
``interesting'' dimensions may only be obtained by taking closed subsets
of $G$ of infinite index.

\subsection{Hausdorff dimension of partition
  subgroups of $G$}\label{sec:hdim-partition}\

The partitions subgroups of $G$ (of infinite index) are of the form
$P_\mu$ for a partition diagram 
$\mu=\{(i,j)\in\N^2~|~1\leq i<j\}\subset\N^2$, 
subject to $(i,j),(j,k)\in\mu\Rightarrow(i,k)\in\mu$; or a partition
$\mu=(\mu_2,\mu_3,\dots)$, where
$\mu_j=\max\{i~|~(i,j)\in\mu\}$ for all $j\geq2$ as defined in
  Section~\ref{sec:partition}.
Hence
$P_\mu=\{x\in G~|~x_{ij}=0\;,\;\forall\;(i,j)\notin\mu\}$, which gives
\begin{equation}\label{eqn:dim-p}
\frac{\log|P_\mu N_n:N_n|}{\log|G_n|}=\frac{2|\mu|_n}{n(n-1)}
\end{equation}
where $|\mu|_n$ denotes the cardinality of the subset of the squares
  of $\mu$ up to, and including, the $n$th column: 
$$|\mu|_n=|\{(i,j)\in\mu|~1\leq i<j\leq n\}|
\qbox{for any integer $n\geq2$.}$$

\smallskip
With this notation, we can state and prove the main result in this
section. 

\begin{thm}\label{thm:hdim-partition}
Let $\alpha\in[0,1]$. Then there exists a partition subgroup
\;$P_\mu$\; for which\;
 $\dim(P_\mu)=\alpha$. 
In particular $\spec(G)=[0,1]$.

Moreover, for all $d\geq1$,
$$\dim(\gamma_d(G))=\dim(G^{(d)})=1$$
for the subgroups of $G$ in the lower central and derived series of
$G$. 
\end{thm}

\begin{proof}
The existence of a partition subgroup for $\alpha\in\{0,1\}$ is
clear. 
Suppose $\alpha\in(0,1)$ and let $(a_n)_{n\geq2}\subset\Q$ be the
sequence of rational numbers defined as follows: 
$$a_n=\frac{2b_n}{n(n-1)}\qbox{where}
b_n=\lfloor\frac{\alpha n(n-1)}2\rfloor
\qbox{for all $n\in\N$,}$$
where $\lfloor x\rfloor\leq x<\lfloor x\rfloor+1$ denotes the integer
part of any $x\in\R$. 

We claim that $(a_n)_{n\geq2}$ converges to $\alpha$.
The inequalities $\dst b_n\leq\frac{\alpha n(n-1)}2<b_n+1$ imply that
$$a_n=\frac{2b_n}{n(n-1)}\leq\alpha<
\frac{2(b_n+1)}{n(n-1)}=a_n+\frac2{n(n-1)}.$$
Let $\varepsilon>0$. We want to show that there exists $N\in\N$ such
that $|\alpha-a_n|<\varepsilon$ for all $n\geq N$. 
Define $N$ as being
the least positive integer such that
$\varepsilon\geq\dst\frac2{N(N-1)}$.
Then for any $n\geq N$, we have
$$|\alpha-a_n|<\frac2{n(n-1)}\leq\varepsilon
\qbox{which proves that}
\lim_{n\to\infty}a_n=\alpha.$$

For $n\geq3$, consider 
\begin{align*}
b_n-b_{n-1}&=\lfloor\frac{\alpha n(n-1)}2\rfloor-
\lfloor\frac{\alpha(n-1)(n-2)}2\rfloor\\
&=\lfloor\frac{\alpha(n-1)(n-2)}2+\alpha(n-1)\rfloor
-\lfloor\frac{\alpha(n-1)(n-2)}2\rfloor
\end{align*}
where $0<\alpha(n-1)$, so that the difference is nonnegative for all
$n\geq2$.
More precisely, from the inclusion
$b_n-b_{n-1}\in\big(\alpha(n-1)-1\;,\;\alpha(n-1)+1\big)$, which
contains a unique integer, we see that whenever
the difference is positive, then $b_n-b_{n-1}$ is of the order of 
$\alpha(n-1)$. 

Now, define the partition $\mu=(\mu_2,\mu_3,\dots)\subset\N$ as
follows: $\mu_2=b_2$ and then $\mu_n=b_n-b_{n-1}$
for all $n\geq 3$. Because the
differences $b_n-b_{n-1}$ are nonnegative integers less than $n$, the
sequence $\mu$ defines a partition subgroup $P_\mu$. That is, $P_\mu$ is
the subgroup of $G$ whose nonzero squares in column 
$n$ are the top $b_n-b_{n-1}$ ones (possibly none).  
For $n\geq2$, we have 
$$|\mu|_n=b_2+\sum_{3\leq j\leq n}b_j-b_{j-1}=b_n\qbox{so that}
|P_\mu N_n/N_n|=q^{b_n}.$$
It follows that
$$\dim(P_\mu)=\lim_{n\to\infty}\frac{\log|P_\mu N_n:N_n|}{\log|G_n|}=
\lim_{n\to\infty}\underbrace{\frac{2b_n}{n(n-1)}}_{=a_n}=\alpha$$
proving the first part of the theorem.

To prove the second part of the statement, recall that
$$\gamma_{d}(G)=P_{(0^{d-1},1,2,3,\dots)}\qbox{and}
G^{(d)}=P_{(0^{(-1+2^{d-1})},1,2,3,\dots)}\qbox{for all $d\geq2$.}$$
So 
\begin{align*}
|\gamma_d(G)N_n/N_n|&=|G_{n+1-d}|=q^{\frac{(n+1-d)(n-d)}2}\qbox{and}\\
|G^{(d)}N_n/N_n|&=|G_{n+1-2^{d-1}}|=q^{\frac{(n+1-2^{d-1})(n-2^{d-1})}2}.
\end{align*}
It follows that 
$$\dim(\gamma_d(G))=
\lim_{n\to\infty}\frac{\log|\gamma_d(G)N_n:N_n|}{\log|G_n|}
=\lim_{n\to\infty}\frac{(n+1-d)(n-d)}{n(n-1)}=1~,$$
and similarly
$$\dim(G^{(d)})=\lim_{n\to\infty}
\frac{(n+1-2^{d-1})(n-2^{d-1})}{n(n-1)}=1
\qbox{for all $d\geq2$.}$$
\end{proof}

\begin{remark}\label{rem:hdim}
From the proof of the theorem, we see that one could try to modify the
definition of $\mu$ in order to obtain a normal partition subgroup with
prescribed Hausdorff dimension. This ``tweaking'' consists in shifting
each ``bulging'' square arising whenever $b_n-b_{n-1}>b_{n+1}-b_n$ by a
finite number of columns to the right until it reaches the next
``landing''. The next examples may shed some light on this.
\end{remark}

\begin{example}
\begin{enumerate}
\item 
Let $\alpha=\dst\frac1{\pi}\approx0.318309886$. The sequence
$(a_n)_{n\in\N}$ obtained by the method of the proof of
Theorem~\ref{thm:hdim-partition} gives the integers
$\dst b_n=\lfloor\frac{n(n-1)}{2\pi}\rfloor$\;and $\mu_n=b_n-b_{n-1}$,
with $\mu_2=b_2$. We calculate
$$\begin{array}{r|r|r|r|r|r|r|r|r|r|r|r|r|r|r|r|r|r|r|r|r}
n&2&3&4&5&6&7&8&9&10&11&12&13&14&15&16&17&18&19&20\\
\hline
\mu_n&0&0&1&2&1&2&2&3&3&3&4&3&4&5&5&5&5&6&6
\end{array}~.$$

That is, the subgroup
$$P_\mu=\left(\begin{array}{rrrrrrrrrrrrrrrr}
1&0&0&*&*&*&*&*&*&*&*&*&*&*&\dots\\
&1&0&0&*&0&*&*&*&*&*&*&*&*&\dots\\
&&1&0&0&0&0&0&*&*&*&*&*&*&\dots\\
&&&1&0&0&0&0&0&0&0&*&0&*&\dots\\
&&&&\ddots&0&\dots&0&\dots&0&\dots&0&\dots&&\dots
\end{array}\right)~.$$
So the proportion of squares in $P_\mu$ (relative to the total number of
squares in $G$) up to the fourth column is
$\dst\frac16$, while up to the twentieth column we
have $\dst\frac{60}{190}=\frac6{19}\approx0.315789$ (an error of
$0.080\%$ to $3$ decimal places). 

Expanding on Remark~\ref{rem:hdim}, we can tweak the partition $\mu$ to
get a normal partition $\mu'=(0^2,1^2,2^3,3^4,4^2,5^4,6,6,\dots)$, such
that $P_{\mu'}\triangleleft G$ and $\dim(P_{\mu'})=\dim(P_\mu)=\pi^{-1}$.

One can check that the same values for $\mu_n$ for $n\leq20$ are
obtained with $\alpha=\dst\frac7{22}$ instead of $\pi^{-1}$.

\item Let $\alpha=e^{-3}\approx0.049787$. The corresponding partition
  $\mu$ is
$$\begin{array}{r|r|r|r|r|r|r|r|r|r|r|r|r|r|r|r|r|r|r|r|r}
n&2&3&4&5&6&7&8&9&10&11&12&13&14&15&16&17&18&19&20\\
\hline
\mu_n&0&0&0&0&0&1&0&0&1&0&1&0&1&1&0&1&1&1&1\\
\end{array}~.$$
So the proportion of squares in $P_\mu$ up to the
twentieth column is $\dst\frac9{190}\approx0.047368$ (an error of
$4.858\%$ to $3$ d.p.). Here finding a normal partition subgroup with
Hausdorff dimension $e^{-3}$ seems more difficult. 

\end{enumerate}
\end{example}

Theorem~\ref{thm:hdim-partition} gives one amongst ``many'' partition
diagrams which give partition subgroups of prescribed Hausdorff dimension. 
Closed subgroups of $G$ with a rational Hausdorff dimension can easily
be described as partition subgroups using (maybe more natural) partition
diagrams.

\begin{example}
To get a partition subgroup with Hausdorff dimension
$\alpha=\dst\frac12$, let
$$\mu=\{(i,i+2)~|~i\in\N\}~.\qbox{Then}
P_\mu=\begin{pmatrix}1&0&*&0&*&0&*&\dots\\
&1&0&*&0&*&0&\dots\\
&&1&0&*&0&*&\dots\\
&&&1&0&*&0&\dots\\
&&&&1&0&*&\dots\\
&&&&&1&0&\dots\\
&&&&&&1&\dots\\
&&&&&&&\ddots
\end{pmatrix}~.$$ 
Note that $P_\mu$ is not normal in $G$, and that $P_\mu$ is formed by a
half of the super diagonals of $G$. Therefore 
$$\dim(P_\mu)=\lim_{n\to\infty}
\frac{\log|P_\mu N_n:N_n|}{|G_n|}=\frac12~.$$

\end{example}

\subsection{Hausdorff dimension of finitely determined closed
  subgroups}\label{sec:finite-det}\

Let $H\leq_cG$, where $G=G(q)$ be determined by a finite total number of
rows, columns and super diagonals. That is, there exists a finite number 
$R_{i_1},\dots,R_{i_r}$ of rows and $D_{j_1},\dots,D_{j_d}$ of super
diagonals such that for any $x\in H$ any coefficient $x_{ij}$ of $x$ is
given as a function of certain coefficients in the rows and diagonals
$R_{i_l},D_{j_m}$ above. We call such a subgroup {\em finitely
  determined}. 
For instance, any $d$-periodic subgroup is finitely determined by its
first $d$ rows: 
if $H=H[d]$, then $x=(x_{ij})$ is subject to the constraints
$x_{ij}=x_{\bar i\bar j}$, where $\bar i\equiv i\pmod d$ and 
$\bar j-\bar i=j-i$. Similarly, a string subgroup which embeds into a
direct product of subgroups isormorphic to
$\dst\prod_{1<n_i<N}G_{n_i}$, for some upper bound $N$ on the size of
the diagonal blocks, is determined by its first $N$ super diagonals. 

In the examples seen above, we observed that the subgroup
$\sigma(\CN)$ isomorphic to the Nottingham group is determined by its
first row (cf. paragraph preceeding Proposition~\ref{prop:cn}), 
while the free product $\varphi(F)$, as a subgroup of $G(p)$, is
determined by its first $2$ rows, by
Proposition~\ref{prop:f-periodic}. 
Recall from Section~\ref{sec:free} that
$$\varphi(F)=\varphi(F)[2]=\langle s=1+\sum_{n\in\N}e_{2n-1,2n}~,
~t=1+\sum_{n\in\N}e_{2n,2n+1}\rangle\cong C_p*C_p~.$$ 
More generally, for a $d$-periodic subgroup $H$, the subgroup
$H_n=HN_n/N_n$ of $G_n$ is the subgroup whose coefficients in the
first $d$ rows can be chosen freely, and these uniquely determine the
remaining ones in the bottom $(n-d)$ rows. 
$$\begin{pmatrix}1&a_{12}&\dots&\dots&\dots&a_{1n}&|&a_{1,n+1}&\dots\\
&\ddots&\dots&\dots&\dots&\vdots&|&\dots&\dots\\
&&1&a_{d,d+1}&\dots&a_{d,n}&|&a_{d,n+1}&\dots\\
&&&\ddots&\bullet\dots&\bullet&|&\bullet_{i,n+1}&\dots\\
&&&&1&\bullet&|&\bullet_{n-1,n+1}&\dots\\
  &&&&&1&|&\bullet_{n,n+1}&\dots\\
\hline
&&&&&&|&\hbox{\Huge*}\end{pmatrix}$$
where $\bullet$ denote the coefficients that are determined by the freely
chosen chosen $a_{ij}$.

Although not finitely generated, nor finitely presented in general,
finitely determined subgroups of $G$ are ``small'' in $G$, in
the following sense.

\begin{lemma}\label{lem:fin-det}
Suppose that $H\leq_cG$ is finitely determined. 
Then $\dim(H)=0$.
\end{lemma}

Note that a closed subgroup of $G$ which
is determined by a finite number of columns is also determined by a
finite number of rows, so that the lemma applies to this
class of subgroups too.

\begin{proof}
We prove the claim for $H\leq_cG$ determined by a finite number of rows
$i_1<\dots<i_d$. Let $H_n=HN_n/N_n$ and $K_d$ the subgroup of $G$ formed
by all the $d$-periodic elements. So $K_d$ is determined by its first
$d$-rows. For $n>d$ we calculate
\begin{align*}
\log|H_n|&\leq\log|K_dN_n/N_n|=\log|G_n|-\log|G_{n-d}|\\
&=\frac{n(n-1)}2-\frac{(n-d)(n-d-1)}2=\frac{d(2n-d-1)}2~.
\end{align*}
It follows that
$$\dim(H)\leq\lim_{n\to\infty}\frac{d(2n-d-1)}{n(n-1)}=0
\qbox{and so}\dim(H)=0.$$
Similarly, if $H$ is determined by $d$ super diagonals, then for $n>d$,
we have (counting the coefficients in the successive super diagonals up
to the $d$th one)
$$\log|H_n|\leq(n-1)+(n-2)+\dots+(n-d)=\frac{d(2n-d-1)}2$$
and we conclude as above.
The lemma follows.
\end{proof}

Here are immediate consequences of this observation.

\begin{cor}\label{cor:fin-dim}
The following hold.
\begin{enumerate}
\item The subgroup $\sigma(\CN)$ of $G$ isomorphic to the Nottingham
  group has Hausdorff dimension $0$.
\item If $q=p$, the subgroup $\varphi(F)\cong C_p*C_p$ of $G$ has
  Hausdorff dimension $0$.
\item If $H$ is a $d$-periodic subgroup of $G$ for some positive integer
  $d$, then $\dim(H)=0$.
\item If $H$ is a string subgroup of $G$ of the form
  $H\leq\dst\prod_{1<n_i\leq N}G_{n_i}$ for some integer $N\geq2$, then
 $H$ has Hausdorff dimension $0$.
\end{enumerate}
\end{cor}

\subsection{Hausdorff dimension and field
  extensions}\label{sec:hdim-field-ext}\ 

We use the same notation as in Section~\ref{sec:field-ext}.
Let $H=\alpha_f(G(q))\leq_cG(p)$ and $q=p^f$ for some $f\in\N$. 
Given $n\in\N$ write $n=rf+s$ with $0\leq s<f$
and $r\geq0$. Write $\log=\log_p$. We calculate

\begin{align*}
\frac{\log|G_r(q)|}{\log|G_n(p)|}&\leq
\frac{\log|HN_n(p):N_n(p)|}{\log|G_n(p)|}\leq
\frac{\log|G_{r+1}(q)|}{\log|G_n(p)|}\\
\frac{\log|q^{\frac{r(r-1)}2}|}{\log|p^{\frac{n(n-1)}2}|}&
\leq\frac{\log|HN_n(p):N_n(p)|}{\log|G_n(p)|}\leq
\frac{\log|q^{\frac{r(r+1)}2}|}{\log|p^{\frac{n(n-1)}2}|}\\
\frac{fr(r-1)}{(rf+s)(rf+s-1)}&\leq
\frac{\log|HN_n(p):N_n(p)|}{\log|G_n(p)|}\leq
\frac{fr(r+1)}{(rf+s)(rf+s-1)}~.
\end{align*}
Left- and right hand side terms both converge to $\dst\frac1f$ as
$n\to\infty$, i.e. as $r\to\infty$. Therefore
$\dst\dim\big(\alpha_f(G(q))\big)=\frac1f$.

\smallskip
Let $K=\beta_f(G(p))\leq_c G(q)$ and $q=p^f$. By definition of
$\beta_f$, for each $n\geq2$, the index of $KN_n(q)/N_n(q)\cong G_n(p)$
in $G_n(q)$ is equal to
$\frac{p^{\frac{n(n-1)}2}}{(p^f)^{\frac{n(n-1)}2}}$, so that 
$$\lim_{n\to\infty}\frac{\log|KN_n(q):N_n(q)|}{\log|G_n(q)|}=\frac1f~.$$

\subsection{Hausdorff spectrum of $G(\Z_p)$}\label{sec:hdim-p-adic}\

To calculate the Haudorff dimension of closed subgroups of $G(\Z_p)$, we
consider the filtration given in Equation~(\ref{eqn:p-adic-base}):
$$V_{n}=\left\{\left(\begin{array}{r|r}\dst
  G_{n}(p^n\Z_p)&\hbox{\Large*}\\
\hline
\hbox{\Large0}&\hbox{\Large*}\end{array}\right)~|~*\in\Z_p\right\}$$
with factor groups
$\dst G/V_{n}\cong G_{n}(\Z/p^n)$ for all $n\in\N$. 

Theorem~\ref{thm:hdim-partition} proves that the dimension spectrum of
$G(\Z_p)$ is the whole interval $[0,1]$, which can be attained using
solely partition subgroups of $G(\Z_p)$. 

\begin{cor}\label{cor:hdim-partition}
Let $\alpha\in[0,1]$. Then there exists a partition subgroup
$P_\mu$ for which $\dim(P_\mu)=\alpha$. 
In particular $\spec(G(\Z_p))=[0,1]$.

For all proper ideals $I\subset\Z_p$ and any ideal partition
subgroup $P_\mu(I)$ as in Proposition~\ref{prop:p-adic-part}, we have
$\dim(P_\mu(I))=0$. 

Moreover, for all $d\geq1$,
$$\dim(\gamma_d(G(\Z_p)))=\dim(G(\Z_p)^{(d)})=1$$
for the subgroups of $G(\Z_p)$ in the lower central and derived series
of $G(\Z_p)$. 
\end{cor}


\appendix

\section{Automorphisms of the finite groups $G_n(q)$}\label{app:aut-gn}

Let $G=G_n(q)$ where $q=p^f$ for some $f,n\in\N$.
So $G_n(q)$ is generated by all the matrices of the form
$1+a_ie_{rs}$, where $a_1,\dots,a_f$ form a basis of $\F_q$ as
$\F_p$-vector space and $1\leq r<s\leq n$.
A minimal set of generators is formed by all such elements of the form
$1+a_ie_{r,r+1}$. In \cite{weir}, A. Weir determines the group of
automorphisms of $G_n(q)$ and describes the maps by their action on the
elements $1+a_ie_{r,r+1}$ in a minimal set of generators.

\begin{thm}\cite[Theorem~8]{weir}\label{thm:weir-aut}
The group $\Aut(G_n(q))$ of automorphisms of $G_n(q)$ is generated by
the subgroups
$\langle\tau\rangle,~\mathcal L,~\mathcal D~,\mathcal I$ and $\mathcal P$, where 
\begin{align*}
\langle\tau\rangle(1+e_{ij})&
=1+e_{n+1-j,n+1-i}\qbox{for all}1\leq i<j\leq n~;\\
\mathcal L&
=\langle\nu~:~1+a_ie_{r,r+1}\mapsto1+\nu(a_i)e_{r,r+1}\;
a_i\in\F_q~,\;1\leq r<n\;,\;\nu\in\Aut(\F_q)\rangle~;\\
\mathcal D&
=\;\hbox{conjugation by non scalar diagonal matrices;}\\
\mathcal I&
=\;\hbox{conjugation by elements of $G_n(q)$;}\\
\mathcal P&
=\mathcal Z\times\mathcal U\qbox{where}\\
\mathcal Z&=
\left\langle
\begin{array}{rl}
\tau^i_r~:&~1+a_ie_{s,s+1}\mapsto1+a_ie_{s,s+1}+\delta_{r,s}
b_ie_{1,n}~,\\
&~b_i\in\F_q~,~1\leq i\leq f~,~1\leq s<n~,2\leq r\leq n-2\end{array}
\right\rangle\qbox{and}\\
\mathcal U&
=\langle1+ae_{12}\mapsto1+ae_{12}+abe_{2n}\;,\;
1+ae_{n-1,n}\mapsto1+ae_{n-1,n}+abe_{1,n-1}\rangle~.
\end{align*}
In particular, $\langle\tau\rangle$ has order $2$, $\mathcal L$ is
cyclic of order $f$, $\mathcal D$ has order $(q-1)^{n-1}$, $\mathcal I$
has order $|G_n(q)/Z(G_n(q))|=q^{\frac{n^2-n-2}2}$, $\mathcal P$ is
elementary abelian of order $q^{f(n-3)+2}$. 
\end{thm}
The elements of $\mathcal Z$ are called {\em central} automorphisms,
because they induce the identity on $G/Z(G)$, and those of $\mathcal U$
are called {\em extremal} automorphisms. The maps
$\tau^i_1,\tau^i_{n-1}$ are inner automorphisms, and so need not be
added in $\mathcal Z$.
The function $\tau$ corresponds to the symmetry of the Dynkin diagram of
type $A_{n-1}$, flipping the squares about the antidiagonal. 

We now prove the technical lemma which we used in the proof of
Proposition~\ref{prop:field-ext-norm}. We use the same notation as in
the proposition. In particular, $\alpha_f~:~G(q)\to G(p)$, induced by
regarding $\F_q$ as an $f$-dimensional $\F_p$-vector space, has image
$H$, and $\beta_f~:~G(p)\to G(q)$, induced by the inclusion of the
coefficients $\F_p\subseteq\F_q$, has image $K$.

\begin{lemma}\label{lem:centr-h}
Let $H_m=HN_m(p)/N_m(p)$ and $L_m=\alpha_f(K)N_m(p)/N_m(p)$ for any
$m\in\N$. Then 
$C_{G_{nf}(p)}(L_{nf})=C_{G_{nf}(p)}(H_{nf})=
\langle1+e_{ij}~|~1\leq i\leq f~,~n-f<j\leq n\rangle$, is the
subgroup formed by all the matrices whose only nonzero nondiagonal
squares lie in the upper right $f\times f$ corner.
\end{lemma}

\begin{proof}
Let $X=\langle1+e_{ij}~|~1\leq i\leq f~,~n-f<j\leq n\rangle$. 
Clearly the elements of $X$ commute with any element
of $H_{nf}$ because the first and last diagonal $f\times f$ blocks of
any element in $H_{nf}$ is the identity $f\times f$ matrix. 

To show that these are exactly the elements which centralise $H_n$, it
is enough to see that no other element centralises an element of
$L_{nf}$. Let
$y=\alpha_f(1+e_{u+1,v+1})=1+\dst\sum_{1\leq k\leq f}e_{uf+k,vf+k}$
for $0\leq u<v< n$ and suppose that 
$\dst x=1+\sum_{1\leq i<j\leq n}a_{ij}e_{ij}\in C_{G_{nf}(p)}(L_{nf})$. 
Put $x^{-1}=1+\dst\sum_{1\leq i<j\leq n}b_{ij}e_{ij}$. We calculate
\begin{equation}\label{eq:conjug}
\ls xy=y+\underbrace{\sum_{1\leq k\leq f}\bigg(\sum_{vf+k<j}\bigg(
b_{vf+k,j}e_{uf+k,j}+\sum_{i<uf+k}a_{i,uf+k}b_{vf+k,j}e_{ij}\bigg)\bigg)}_{(*)}
\end{equation}
and solve the equation $(*)=0$. Note that all the indices $(i,j)$
appearing in $(*)$ are distinct. Therefore $b_{vf+k,j}=0$ for all
$j>vf+k$, all $1\leq k\leq f$ and all $1\leq v<n$. By definition of these
coefficients, we must then also have $a_{ij}=0$ for all $f<i<j\leq n$.
Now, take 
$y=\alpha_f(1+e_{12})=
\begin{pmatrix}I_f&I_f&0\\&I_f&0\\&&I_{(n-2)f}\end{pmatrix}$
and suppose that
$x=\begin{pmatrix}A&B&C\\&I_f&0\\&&I_{(n-2)f}\end{pmatrix}$ centralises
$y$. We calculate 
$\ls xy=\begin{pmatrix}I_f&A&0\\&I_f&0\\&&I_{(n-2)f}\end{pmatrix}$, which
gives $A=I_f$. Similarly, for $y=\alpha_f(1+e_{23})$ we obtain that $B$
is the zero $f\times f$ matrix.
Inductively on $\alpha_f(1+e_{r,r+1})$ each of the successive 
$f\times f$ blocks of the $f\times(n-2)f$ matrix $C$ must be zero,
except the last one in which the squares can take any value. This proves
the lemma.
\end{proof}


\section{Fractional dimension for profinite groups}\label{app:frac-dim}

We review the concept of fractional dimension for profinite groups, as
introduced in \cite{aber}, and refer to \cite{falconer} for the
concepts used.

First some background on measure theory. Let $X$ be a set and write
$\mathcal P(X)$ for the power set of $X$. A non-empty subset  
$\mathcal S\subseteq\mathcal P(X)$ is a $\sigma$-field if $\mathcal S$
is closed under taking complements and countable unions. 
The {\em Borel sets} of $X$ are the sets belonging to the $\sigma$-field
generated by the closed subsets of $X$.
Given a
$\sigma$-field, one can show that
$$\liminf{j\to\infty}E_j=
\bigcup_{k\geq1}\bigcap_{j\geq k}E_j\in\mathcal S
\qbox{and}\limsup_{j\to\infty}E_j=
\bigcap_{k\geq1}\bigcup_{j\geq k}E_j\in\mathcal S~.$$
The former set is formed by all elements which are in all but a finite
number of $E_j$, while the latter set is formed by all elements which
belong to infinitely many $E_j$.

A {\em measure} defined on a $\sigma$-field $\mathcal S$ is a function
$\mu~:~\mathcal S\to[0,\infty]$
such that $\mu(\emptyset)=0$ and $\mu(\cup_jE_j)=\sum_j\mu(E_j)$ for
any countable collection of disjoint sets $E_j$. We call $\mu$ an {\em
  outer measure} if $\mathcal S=\mathcal P(X)$, and we call it a {\em
  probability measure} if $\mu(X)=1$. 

Suppose that $X$ is a profinite group $X=\dst\prlim{n\in\N}X_n$ with
projection maps $\theta_n~:~X\to X_n$ onto the finite quotients and maps
$\pi_{n,m}~:~X_n\to X_m$ such that $\theta_m=\pi_{n,m}\theta_n$ for all
$m\leq n$ (with $\pi_{n,n}=\id_{X_n}$). We consider the {\em standard
  basis} 
$$\mathcal B=\{\theta_n^{-1}(x_n)~|~x_n\in X_n~,~n\in\N\}$$
for the topology on $X$. For short, let $N_n=\ker(\theta_n)$ for all
$n\in\N$. It is well-known in measure theory that the so-called Haar
measure is the unique probability measure on $X$ which is
$X$-invariant. This measure satisfies 
$$\mu(\theta_n^{-1}(x_n))=\mu(N_n)=|X_n|^{-1}\qbox{and}
\mu(\{x\})=0\;\forall\;x\in X.$$
This second property means that $\mu$ is {\em non-atomic}. 

Now define $\Delta_\mu~:~\mathcal P(X)\to\R$ as follows.
For all $\delta,\gamma\in\R$ with $\delta>0$ define for 
$M\in\mathcal P(X)$,
$$l_{\mu,\theta}^\gamma(M)=
\inf_{\mathcal C}\sum_{B\in\mathcal C}(\mu(B))^\gamma$$
where $\mathcal C$ is a cover of $M$ by balls $B\in\mathcal B$ such that
$\mu(B)<\theta$. Hence let
$$l_\mu^\gamma=\lim_{\theta\to0}l_{\mu,\theta}^\gamma~.$$
One can show that there exists a unique real number $\Delta_\mu(M)$
satisfying
$$l_\mu^\gamma(M)=\infty\qbox{for all $\gamma<\Delta_\mu(M)$, and}
l_\mu^\gamma(M)=0\qbox{for all $\gamma>\Delta_\mu(M)$~.}$$
Thus $\Delta_\mu(M)$ is the {\em Hausdorff dimension of $M$}, or {\em
  Billingsley dimension of $M$} (\cite{aber}).

Abercrombie proves the following.
\begin{prop}
$\Delta_\mu$ is a positive increasing set function and for each
  $\gamma\in\R$, the function $l_\mu^\gamma$ is an outer measure on
  $X$. Furthermore,
\begin{itemize}
\item[(i)] For all $(M_n)_{n\in\N}\subset\mathcal P(X)^{\N}$, then
$\Delta_\mu(\cup_nM_n)=\sup n\big(\Delta_\mu(M_n)\big)$.
\item[(ii)] If $M\in\mathcal B$ such that $\Delta_\mu(M)>0$, then
  $\Delta_\mu(M)=1$.
\item[(iii)] Let $Y$ be a closed subset of $X$, i.e. 
$Y=\dst\bigcap_{n\in\N}YN_n/N_n$. Then
$$\Delta_\mu(Y)\geq\liminf{n\to\infty}\frac{\log|YN_n:N_n|}{\log|X_n|}
\qbox{with}
\Delta_\mu(Y)=\lim_{n\to\infty}\frac{\log|YN_n:N_n|}{\log|X_n|}$$
whenever the limit exists.
\end{itemize}
\end{prop}

\noindent{\bf Acknowledgements.}\;
We wish to thank the referee for a careful proofreading, correcting a
few errors, and for helpful suggestions.

\end{document}